\documentclass[11pt,letter]{article}
\usepackage{amsfonts}

\def\up#1{\raise 1ex\hbox{\small #1}}
\usepackage[matrix,arrow]{xy}
\xyoption{all}

\begin{document}
\bigskip

\centerline{\bf Dualit\'e de Van den Bergh}

\medskip

\centerline{ \bf et Structure de Batalin-Vilkoviski\v{\i} sur les alg\`ebres de Calabi-Yau}
\medskip

\bigskip

\centerline{\bf Thierry Lambre
\footnote{Laboratoire de Math\'ematiques,
UMR 6620 du CNRS,
Universit\'e B. Pascal,
63177 Aubi\`ere Cedex.
\hbox{Courriel : thierry.lambre@math.univ-bpclermont.fr}}
}
\bigskip

{\bf R\'esum\'e : } 
Nous montrons que la notion de  calcul de Tamarkin-Tsygan \`a dualit\'e permet de construire des structures de Batalin-Vilkoviski\v{\i} dans un cadre g\'en\'eral. Nous montrons que  la dualit\'e de Van den Bergh des alg\`ebres est  un calcul de Tamarkin-Tsygan \`a dualité. Ceci  permet notamment de retrouver la structure BV des alg\`ebres de Calabi-Yau mise en \'evidence par V. Ginzburg.

\bigskip

{\bf Summary : }
The abstract notion of Tamarkin-Tsygan calculus with duality gives Batalin-Vilkoviski\v{\i} structures in a general setting. We apply this technique to the case of Van den Bergh duality for algebras to prove that  Calabi-Yau algebras are BV-algebras.
\bigskip

{\bf Classification AMS :} 16 E 40, 20 J 06, 55 U 30.
\bigskip

\bigskip
\bigskip
\bigskip

{\bf Introduction.}
\medskip

V. Ginzburg a montr\'e r\'ecemment  que les alg\`ebres de Calabi-Yau sont des alg\`ebres de Batalin-Vilkoviski\v{\i} (en abr\'eg\'e, BV-alg\`ebres).

\medskip
{\bf Th\'eor\`eme} ([Gi, 3.4.3]).
{\sl Soit $A$ une alg\`ebre de Calabi-Yau et soit $(H^*(A,A),\cup,[\ ,\ ])$ l'alg\`ebre de cohomologie de Hochschild de $A$, munie de sa structure d'alg\`ebre de Gerstenhaber. Il exis\-te un g\'en\'erateur $\Delta$ du crochet de Gerstenhaber $[\ ,\ ]$, c'est-\`a-dire qu'il existe une application
$\Delta: H^*(A,A)\to H^{* -1}(A,A)$ satisfaisant \`a la relation 
$$ [\alpha,\beta]=
\Delta(\alpha\cup\beta)-(-1)^p\alpha\cup\Delta(\beta)-\Delta (\alpha)\cup\beta$$ pour tout $\alpha\in H^p(A,A)$ et $\beta\in H^q(A,A)$.}

\medskip
Afin de commenter ce r\'esultat, introduisons quelques notations. 
\smallskip
 
Notons 
$B:H_*(A,A)\to H_{*+1}(A,A)$ le bord de Connes  en {\sl homologie} de Hochschild ([C], [L]).
Soit  $d$ la dimension cohomologique de l'alg\`ebre de Calabi- Yau et soit  $$VdB:H^*(A,A)\to H_{d-*}(A,A)$$ l'isomorphisme de dualit\'e de Van der Bergh ([VdB]).
 
\smallskip
 
Soient  $\alpha\in H^p(A,A)$, $\beta\in H^q(A,A)$ et $z\in H_r(A,A)$. 
La contraction 
$$\iota_\alpha: H_r(A,A)\to H_{r-p}(A,A)$$ est d\'efinie par 
$\iota_\alpha (z)=z\cap \alpha$ 
(le symbole $\cap$ d\'esigne ici le 
{\sl cap}-produit, voir paragraphe 2).
\smallskip

Les deux ingr\'edients essentiels de la d\'emonstration du th\'eor\`eme de V. Ginzburg sont les suivants.

 1) La contraction satisfait l'identit\'e remarquable suivante de Tamarkin-Tsygan ([T-T]).
$$\iota_{[\alpha,\beta]}=
[[B,\iota_{\alpha} ]_{gr},\iota_{\beta} ]_{gr},\leqno (TT)$$
où dans le membre de droite de cette expression, les crochets $[\ ,\ ]_{gr}$ d\'esignent des commutateurs gradu\'es.

\smallskip

2) L'inverse $D=(VdB)^{-1}$ de l'isomorphisme de dualit\'e de Van den Bergh satisfait l'identit\'e remarquable de Ginzburg ([Gi], 3.4.3).
$$D(z\cap\alpha)= D(z)\cup\alpha.\leqno(G)$$

\smallskip
 
Soit $d$ la dimension de l'alg\`ebre de Calabi-Yau. A partir des identit\'es remarquables $(TT)$ et $(G)$, il est facile de v\'erifier qu'en posant $$\Delta=(-1)^dDBD^{-1},$$ on a la relation de Batalin-Vilkoviski\v{\i}
$$[\alpha,\beta]=\Delta(\alpha\cup\beta)-(-1)^p\alpha\cup\Delta(\beta)-\Delta (\alpha)\cup\beta.\leqno (BV)$$
\medskip

Autrement dit le th\'eor\`eme de V. Ginzburg affirme que le conjugu\'e du bord de Connes en homologie de Hochschild par l'inverse de l'isomorphisme de dualit\'e de Van der Bergh est un g\'en\'erateur du crochet de Gerstenhaber de l'alg\`ebre de cohomologie de Hochschild $H^*(A,A)$ de l'alg\`ebre de Cabali-Yau $A$.  
\smallskip

Nous \'enonçons et employons dans ce texte  une g\'en\'eralisation de ce ph\'enom\`ene  pour les calculs de Tamarkin-Tsygan \`a dualit\'e, objets satisfaisant dans un certain cadre de g\'en\'eralit\'e les relations $(TT)$ et $(G)$.
\medskip

{\bf Th\'eor\`eme 1.6.}
{\sl Soit $(H^*, H_*,\kappa, c)$ un calcul de Tamarkin-Tsygan \`a dualit\'e, de classe fondamentale  $c\in H_d$ (voir d\'efinition 1.3).  Notons $D=(c\cap\ -\ )^{-1}$ l'inverse de l'isomorphisme de dualit\'e. Alors un g\'en\'erateur du crochet de Gerstenhaber de $H^*$ est $\Delta=(-1)^dD\kappa D^{-1}$ et l'alg\`ebre de Gerstenhaber $H^*$ est une $BV$-alg\`ebre.}

\medskip

Pour exploiter ce r\'esultat g\'en\'eral dans le cadre de la dualit\'e de Van den Bergh des alg\`ebres, nous montrons que l'isomorphisme de dua\-li\-t\'e de Van den Bergh s'exprime comme le {\sl cap}-produit par une certaine classe fondamentale canoniquement associ\'ee aux alg\`ebres consid\'er\'ees.
\medskip

{\bf Th\'eor\`eme 4.2.}
{\sl Soit  $A$ une alg\`ebre \`a dualit\'e de Van den Bergh, de module dualisant ${\cal D}=H^d(A, A^e)$, de classe fondamentale 
$c\in H_d(A,{\cal D})$. Alors, pour tout $A^e$-module $M$ et pour tout entier $p\geq 0$, le cap-produit 
$$c\cap- :H^p(A, M)\to H_{d-p}(A, {\cal D}\otimes_AM)$$
est un isomorphisme.}
\medskip

Les r\'esultats 1.6 et 4.2 permettent de fournir une d\'emonstration du th\'eor\`eme de V. Ginzburg.
\medskip

Dans diverses situations analogues, une dualit\'e en terme de {\sl cap}-produit est bien connue. C'est le cas notamment en cohomologie des groupes ([B-E]) mais aussi dans un cadre de g\'eom\'etrie de Poisson (voir par exemple [H] et [X]).
Dans ces diff\'erents contextes (groupes, alg\`ebres ou g\'eom\'etrie de Poisson), le {\sl  cap}-produit par une certaine classe fondamentale est un isomorphisme.

En ce sens, les alg\`ebres de Calabi-Yau appara\^{\i}ssent  comme  l'analogue alg\'ebrique des groupes \`a dualit\'e de Poincar\'e {\sl orientables}. Pour se convaincre de la pertinence de cette affirmation, il peut être utile de s'aider du dictionnaire analogique suivant.

\bigskip

\centerline{Groupes \`a dualit\'e 
\ \ \ \ \ \ \ \ \ \
 Alg\`ebres \`a dualit\'e} 
\centerline{${\mathbb{Z}}
\ \ \ \ \ \ \ \ \ \
A$}
\centerline{$\ \ \ \ \ \ \ \ \ \ {\mathbb{Z}}[G]
\ \ \ \ \ \ \ \ \ \ \ \ 
A^e=A\otimes_k A^{op}$}

\bigskip

\centerline{Cohomologie des groupes
\ \ \ \ \ \ \ \ \ \
Cohomologie des alg\`ebres}
\centerline{
$H^p(G,M)=\hbox{Ext}^p_{\mathbb{Z}[G]}({\mathbb{Z}},M)$
\ \ \ \ \ \ \ \ \ \
$H^p(A, M)=\hbox{Ext}^p_{A^e}(A, M)$}

\bigskip

\centerline{Groupe de type FP
\ \ \ \ \ \ \ \ \ 
Alg\`ebre de type FP}
\bigskip
\centerline {Dimension cohomologique}
\centerline{$d=\hbox{pdim}_{{\mathbb{Z}}[G]}({\mathbb{Z}})$
\ \ \ \ \ \ \ \ \ \
$d=\hbox{pdim}_{A^e}(A)$}

\bigskip

\centerline{Module dualisant}
\centerline{${\cal D}=\hbox {Ext}^d_{{\mathbb{Z}}[G]}({\mathbb{Z}},{\mathbb{Z}}[G])
\ \ \ \ \ \ \ \ \ \
{\cal D}=\hbox{Ext}^d_{A^e}(A, A^e)$} 

\bigskip

\centerline{Classe fondamentale}
\centerline{$H_d(G,{\cal D})\cong 
\hbox{Hom}_{{\mathbb{Z}}[G]}({\cal D},{\cal D})
\ \ \ \ \ \ \ \ \ \
H_d(A,{\cal D})\cong \hbox{Hom}_{A^e}({\cal D},{\cal D})$} 

\bigskip

\centerline{Th\'eor\`eme de dualit\'e}
\centerline{
Bieri-Eckmann (1973)
\ \ \ \ \ \ \ \  Van den Bergh (1998)
}\centerline{
si $G$ de type FP
\ \ \ \ \ \ \ \ \ 
 si $A$ de type FP}
 \centerline{
 si $G$ de dimension cohomologique $d$\ \ \ \ \ \ \ \ \ \ 
 si $A$ de dimension cohomologique $d$}
\centerline{
si $\hbox {Ext}^i_{{\mathbb{Z}}[G]}({\mathbb{Z}},{\mathbb{Z}}[G])=0, \ i\not=d$
\ \ \ \ \ \  
si $\hbox{Ext}^i_{A^e}(A, A^e)=0, \ i\not=0$}
\centerline{si ${\cal D}$ est sans torsion sur ${\mathbb{Z}}$
\ \ \ \ \ 
si ${\cal D}$ est $A^e$-inversible}
\centerline{alors} 
\centerline{le cap-produit par la classe fondamentale est un isomorphisme}
\bigskip
\centerline{Groupe \`a dualit\'e de Poincar\'e orientable 
\ \ \ \ \ \ \ \ \ \ \ 
Alg\`ebre de Calabi-Yau}
\centerline{${\cal D}\cong{\mathbb{Z}}$ comme ${\mathbb{Z}}[G]$-module
\ \ \ \ \ \ \ \ \ 
${\cal D}\cong A$ comme $A^e$-module}
\centerline{ avec action triviale de G 
\ \ \ \ \ \ \ \ \ 
avec action triviale de $A^e$}

\bigskip

Ce texte est organis\'e comme suit.
\medskip

1. Structure BV pour les calculs de Tamarkin-Tsygan \`a dualit\'e.

\medskip

2. Rappels sur les structures multiplicatives en th\'eorie de Hochschild.
\medskip

3. Alg\`ebres de type $FP$.
\medskip

4. Dualit\'e de Van den Bergh.

\medskip

5. Une d\'emonstration du th\'eor\`eme de V. Ginzburg.
\medskip

6. Une d\'emonstration d'un r\'esultat de M. Kontsevich.

\bigskip
\bigskip
\bigskip
\bigskip

{\bf 1. Structure BV pour les calculs de Tamarkin-Tsygan \`a dualit\'e.}
\medskip

{\bf D\'efinition 1.1.}
{\sl  Soit $H^*=\oplus_{p \geq 0}H^p$ un $k$-espace vectoriel gradu\'e. On note ${\cal L}^*$ le $k$-espace vectoriel gradu\'e 
d\'ecal\'e, ${\cal L}^p=H^{p+1}$. Pour $\alpha\in H^p$, on pose $\mid\alpha\mid=p$ et $deg(\alpha)=p-1$. 

On dit que 
$H^*$ est une alg\`ebre de Gerstenhaber s'il existe des op\'erations $\cup$ et $[\ ,\ ]$ telles que :
\smallskip

1) $(H^*,\cup)$ est une alg\`ebre gradu\'ee commutative, c'est-\`a-dire une alg\`ebre dans laquelle on a la relation 
$\alpha\cup\beta=
(-1)^{\mid\alpha\mid\cdot\mid\beta\mid}\beta\cup\alpha$.
\smallskip

2) $({\cal L}^*,[\ ,\ ])$  est une alg\`ebre de Lie gradu\'ee, c'est-\`a-dire qu'on a la relation d'antisym\'etrie gradu\'ee
 $$[\alpha,\beta]=
(-1)^{deg(\alpha)\cdot deg(\beta)}[\beta,\alpha]$$
et l'identit\'e de Jacobi gradu\'ee  
$$\displaylines{
(-1)^{deg(\alpha)deg(\gamma)}[\alpha,[\beta,\gamma]]
+(-1)^{deg(\beta)deg(\alpha)}[\beta,[\gamma,\alpha]]
+(-1)^{deg(\gamma)deg(\beta)}[\gamma,[\alpha,\beta]]=0.
}$$

\smallskip

3) Pour tout $\alpha\in {\cal L}^p$, $[\alpha,\ -\ ]$ est une d\'erivation de degr\'e $deg(\alpha)$ de l'alg\`ebre gradu\'ee commutative $(H^*,\cup)$, c'est-\`a-dire qu'on a la relation
$$[\alpha,\beta\cup\gamma]=
[\alpha,\beta]\cup\gamma+
(-1)^{deg(\alpha)\mid\beta\mid}\beta\cup[\alpha,\gamma].$$
}

\medskip
 
{\bf D\'efinition 1.2.}
{\sl Un calcul de Tamarkin-Tsygan est la donn\'ee d'un triplet $(H^*,H_*,\kappa)$ d'espaces vectoriels gradu\'es satisfaisant aux conditions a), b) et c) ci-dessous.
\smallskip

a) $(H^*,\ \cup,\ [\ ,\ ])$ est une alg\`ebre de Gerstenhaber telle que $k\subset H^0$.
\smallskip

b) $H_*$ est un $(H^*,\ \cup)$-module gradu\'e, c'est-\`a-dire qu'il existe une application $k$-lin\'eaire
 
$$H^p\otimes_kH_r\to H_{r-p}$$
$$\alpha\otimes z\to z\cap \alpha$$

telle que pour $z\in H_r$ et $\alpha\in H^p$, en posant 
$\iota_\alpha(z)=(-1)^{rp}\ z\ \cap\alpha$ , on a la relation 
$\iota_\alpha\circ\iota_\beta=
\iota_{\alpha\cup \beta}$.
\smallskip

c) Il existe une application 
$\kappa:H_*\to H_{*+1}$ telle que $\kappa^2=0$ et telle qu'en posant 
$$L_\alpha=
[\kappa,\iota_\alpha]_{gr}=
\kappa\ \iota_\alpha-(-1)^{\mid\alpha\mid}\iota_\alpha\ \kappa,$$ 
la relation suivante est satisfaite : 

$$[ L_\alpha,\iota_\beta]_{gr}=
\iota_{[\alpha,\beta]}.\leqno(TT)$$
}
\medskip
 
 {\sc Remarque :} Dans un calcul de Tamarkin-Tsygan $(H^*,H_*,\kappa)$, l'espace vectoriel gradu\'e $H_*$ est un 
$({\cal L}^*,[\ ,\ ])$-module de Lie gradu\'e pour l'op\'eration 
$${\cal L}^{p-1}\times H_r\to H_{r-(p-1)}$$
$$\alpha\otimes z\to  L_\alpha(z),$$
c'est-\`a-dire qu'on a la relation 
$L_{[\alpha,\beta]}
=[L_\alpha,L_\beta]_{gr}
=L_\alpha L_\beta-(-1)^{deg(\alpha)deg(\beta)}L_\beta L_\alpha.$

\medskip

{\bf D\'efinition 1.3.}
{\sl Soit $(H^*,H_*,\kappa)$ un calcul de Tamarkin-Tsygan. On dit que ce calcul est un calcul de Tamarkin-Tsygan \`a dualit\'e s'il existe $c\in H_d$ tel que $c\cap 1=c$ et tel que pour tout entier $p$, 
$$c\cap\ -: H^p\to H_{d-p}$$ est un isomorphisme.
\smallskip

Cet isomorphisme est  appel\'e isomorphisme de dualit\'e. L'\'el\'ement $c$ est appel\'e classe fondamentale du calcul de Tamarkin-Tsygan \`a dualit\'e.}
\medskip

{\bf Proposition 1.4.} (formule de Ginzburg, [G], 3.4.3. ($i$))
{\sl Soit $(H^*, H_*,\kappa, c)$ un calcul de Tamarkin-Tsygan \`a dualit\'e, de classe fondamentale $c\in H_d$. 
Soit $D=(c\cap\ -\ )^{-1}$ l'inverse de l'isomorphisme de dualit\'e. Alors pour tout $\alpha\in H^p$ et tout $z\in H_r$, dans $H_{r-p}$, on a l'\'egalit\'e 
$$D(z\cap\alpha)=D(z)\cup\alpha.\leqno (G)$$}
\medskip

D\'emonstration. Soit $z\in H_r$. Par l'isomorphisme  de dualit\'e, il existe $\beta\in H^{d-r}$ tel que 
$z=c\cap\beta $, ce qui s'\'ecrit encore $D(z)=\beta$.
Pour tout $\alpha\in H^p$, on a  donc 
$$\displaylines{
z\cap\alpha=(c\cap\beta)\cap\alpha\hfill\cr
\ \ \ \ \ \ \  =(-1)^{rp}\ \iota_\alpha(c\cap\beta)\hfill\cr
\ \ \ \ \ \ \ \ \  =(-1)^{rp}\ (-1)^{d(d-r)}\ \iota_\alpha\ \iota_\beta\ (c)\hfill\cr
\ \ \ \ \ \ \ \ =(-1)^{rp+d(d-r)}\ \iota_{\alpha\cup\beta}\ (c)\hfill\cr
\ \ \ \ \ \ \ \ =(-1)^{rp+d(d-r)}\ (-1)^{(p+d-r)d}\ c\cap(\alpha\cup\beta)\hfill\cr
\ \ \ \ \ \ \ \  =(-1)^{rp+d(d-r)+(p+d-r)d}\ (-1)^{p(d-r)}c\cap(\beta\cup\alpha)\hfill\cr
\ \ \ \ \ \ \ \ =
c\cap (D(z)\cup\alpha),\hfill\cr
}$$
ce qui s'\'ecrit encore 
$D(z\cap\alpha)=D(z)\cup\alpha.$
\medskip
 
{\bf Lemme 1.5.}
{\sl Soit $(H^*, H_*,\kappa,c)$ un calcul de Tamarkin-Tsygan \`a dualit\'e, de classe fondamentale $c\in H_d$. Soit $D=(c\cap\ -\ )^{-1}$ l'inverse de l'isomorphisme de dualit\'e.
On d\'efinit $\Delta: H^{*}\to H^{*-1}$ par 
$$\Delta= (-1)^d\ D\  \kappa\  D^{-1}.$$
Soient $z\in H_r$, $\alpha\in H^p$ et 
$\beta\in H^q$.
Alors on a la relation 
$$\displaylines{
D(z)\cup[\alpha,\beta]=
(-1)^{(d-r)}\ \Delta\big( D(z)\cup\alpha\cup\beta\big)\ 
-(-1)^{p(r+d+1)+d-r}\ \alpha\cup\Delta\big( D(z)\cup\beta\big)
\hfill\cr
-(-1)^{(d-r)}\ \Delta\big( D(z)\cup\alpha\big)\cup\beta\ 
-(-1)^{pq-1+d-r}\ \Delta \big (D(z)\big)\cup(\alpha\cup\beta).\cr
}$$}
 
\medskip
  
D\'emonstration. D'apr\`es la formule de Ginzburg 1.4, on a 
$$D(z\cap[\alpha,\beta])= D(z)\cup[\alpha,\beta].$$ 

Par ailleurs, 
$z\cap[\alpha,\beta]=(-1)^{(p+q-1)r}
\ \iota_{[\alpha,\beta]}(z).$
On a donc 
$$D\big( z\cap\ [\alpha,\beta]\big)=
\ (-1)^{(p+q-1)r}\ D\big(\iota_{[\alpha,\beta]}(z)\big).$$ 
La relation $(TT)$ de 1.2 s'\'ecrit  
$$\iota_{[\alpha,\beta]}=
[[\kappa,\iota_\alpha]_{gr},\iota_\beta]_{gr}$$

 d'où
$$
D(z\cap[\alpha,\beta])=(-1)^{(p+q-1)r}D(\iota_{[\alpha,\beta]}(z))
=(-1)^{(p+q-1)r}\ (A-B-C-D)$$
avec 

$A=D\ \kappa\ \iota_{\alpha\cup\beta}\ (z),$

$B=(-1)^p\ D\ \iota_\alpha\  \kappa\ \iota_\beta\ (z),$

$C=(-1)^{(p-1)q}\ D\ \iota_\beta\ \kappa\ \iota_\alpha\ (z)$ 

et $D=(-1)^{(p-1)(q+1)}\ D\ \iota_{\alpha\cup\beta}\ \kappa\  (z).$
\smallskip

Calculons successivement ces quatre termes  en utilisant sans  cesse la relation 
$\Delta D=(-1)^d\ D\kappa$ ainsi que  la formule de Ginzburg.
Les calculs conduisent \`a 
$$\displaylines{
A=D\ \kappa\ \iota_{\alpha\cup\beta}(z)\hfill\cr
=(-1)^d\ (-1)^{(p+q)r}\Delta \big( D(z\cap(\alpha\cup\beta)\big))\hfill\cr
=(-1)^{d+(p+q)r}\Delta(D(z)\cup\alpha\cup\beta).\hfill}$$

\smallskip

De mani\`ere analogue
$$\displaylines{
B=(-1)^p\ D\ \big(\iota_\alpha\ (\kappa\ \iota_\beta\ (z))\big)\hfill\cr
=(-1)^p\ (-1)^{(r-q+1)p}\ D(\kappa\  \iota_\beta\  (z)\cap\alpha)\hfill\cr
=(-1)^{p+(r-q+1)p}\ \big( D\ \kappa\ \iota_\beta\ (z)\big)\cap\alpha\hfill\cr
=(-1)^{p+(r-q+1)p}\ (-1)^d\Delta D(\iota_\beta(z))\cup\alpha\hfill\cr
=(-1)^{p+(r-q+1)p+d}\ (-1)^{rq}\Delta D(z\cap\beta)\cup\alpha\hfill\cr
=(-1)^{p+(r-q+1)p+d+rq}\Delta(D(z)\cup\beta)\cup\alpha\hfill\cr
=(-1)^{p+(r-q+1)p+d+rq}\ (-1)^{(d-r+q-1)p}\alpha\cup\Delta(D(z)\cup\beta)\hfill\cr
=(-1)^{p+dp+rq+d}\ \alpha\cup\Delta(D(z)\cup\beta).\hfill\cr
}$$
 
De mani\`ere analogue, on obtient
$$\displaylines{
C= (-1)^{r(p-q)+d}\Delta(D(z)\cup\alpha)\cup \beta.
}$$
  
\smallskip

Enfin, on a 
$$\displaylines{
D=(-1)^{(p-1)(q+1)}\ D\ \iota_{\alpha\cup_\beta}\ \kappa(z)\hfill\cr
=(-1)^{(p-1)(q+1)}\ (-1)^{(r+1)(p+q)}\ D\big(\kappa(z)\cap(\alpha\cup\beta)\big)\hfill\cr
=(-1)^{(p-1)(q+1)+(r+1)(p+q)}\ D(\kappa(z))\ \cup\ (\alpha\cup\beta)\hfill\cr
=(-1)^{(p-1)(q+1)+(r+1)(p+q)}\ (-1)^d\ \Delta\big(D(z)\big)\ \cup (\alpha\cup\beta).\hfill\cr
}$$
 
\smallskip
 
Ces quatre calculs ach\`event la d\'emonstration du lemme 1.5.

\medskip

{\bf Corollaire 1.6.}
{\sl Soit $(H^*, H_*, \kappa,c)$ un calcul de Tamarkin-Tsygan \`a dualit\'e. Alors 
$H^*$ est une BV-alg\`ebre. Plus pr\'ecis\'ement, en supposant $c\in H_d$ et en notant $D=(c\cap\ -\ )^{-1}$ l'inverse de l'isomorphisme de dualit\'e, un  g\'en\'erateur du crochet de Gerstenhaber de $H^*$ est 
$$\Delta=(-1)^d \ D\kappa D^{-1}$$ et la relation
$$[\alpha,\beta]=
\Delta(\alpha\cup\beta)-(-1)^p\alpha\cup\Delta(\beta)-\Delta (\alpha)\cup\beta$$ est satisfaite.}
\medskip

D\'emonstration. On applique la formule 1.5 \`a $z=c$, classe fondamentale du calcul de Tamarkin-Tsygan. Gr\^ace \`a $c\in H_d$, $D(c)=1$ et $\Delta(1)=0$, on aboutit \`a 
$$[\alpha,\beta]=
\Delta(\alpha\cup\beta)-(-1)^p\alpha\cup\Delta(\beta)-\Delta (\alpha)\cup\beta,$$
 ce qui montre  que $\Delta$ est un g\'en\'erateur du crochet de Gerstenhaber de $H^*$.
 
\bigskip

{\bf 2. Rappels sur les structures multiplicatives en th\'eorie de Hochschild.}
 
\medskip

Soient $k$ un corps et $A$ une $k$-alg\`ebre. On pose $A^e=A\otimes_k A^{op}$.
Soit $M$ un $A^e$-module \`a gauche. 
La cohomologie et l'homologie de Hochschild de $A$ \`a valeur dans $M$ sont  res\-pec\-ti\-vement donn\'ees par 
$H^*(A,M)=\hbox{Ext}^*_{A^e}(A,M)$ et 
$H_*(A,M)=\hbox{Tor}_*^{A^e}(A,M)$. On sait qu'en posant 
$C^p(A,M)=Hom_k(A^{\otimes p},M)$, on a 
$H^p(A,M)=H^p(C^*(A,M), b)$ où $b: C^p(A,M)\to C^{p+1}(A,M)$ est le bord de Hochschild d\'efini pour $f\in C^*(A,M)$  par
$$\displaylines{
bf(a_1,\cdots, a_{p+1})=
a_1f(a_2,\cdots,a_{p+1})
-f(a_1a_2,a_3,\cdots,a_{p+1})+\cdots\hfill\cr
+(-1)^{p}f(a_1,\cdots,a_{p-1},a_{p}a_{p+1})
+(-1)^{p+1}f(a_1,\cdots, a_p)a_{p+1}.
}$$

\smallskip
 
De mani\`ere analogue, en posant $C_r(A,M)=M\otimes_kA^{\otimes r}$, on a 
$H_r(A,M)=H_r(C_*(A,M),b)$ où $b:C_r(A,M)\to C_{r-1}(A,M)$ est donn\'e par la formule 
$$\displaylines{
b(m,a_1,\cdots, a_{r})=
(ma_1,a_2,\cdots,a_{r})
-(m,a_1a_2,a_3,\cdots,a_{r})+\cdots\hfill\cr
+(-1)^{r-1}f(m,a_1,\cdots,a_{r-2},a_{r-1}a_{r})
+(-1)^{r}(a_rm, a_1,\cdots, a_{r-1}).
}$$

\smallskip

Rappelons les d\'efinitions des diff\'erents produits pr\'esents.

\medskip
 
{\sl Le cup-produit.} 
 
Il s'agit d'une application $k$-lin\'eaire 
$$\cup: H^p(A,M)\otimes_k H^q(A,N)\to H^{p+q}(A, M\otimes_AN).$$
 
\'Ecrivons $\alpha\in H^p(A,M)$ sous la forme $\alpha=cl(f)$ et de m\^eme  $\beta\in H^q(A,N)$ sous la forme $\beta=cl(g)$  où $f\in C^p(A,M)$ et  
$g\in C^q(A,N)$ sont des $b$-cocycles. On d\'efinit l'\'el\'ement  
$f\cup g$ de $C^{p+q}(A, M\otimes_AN)$ par la formule
$$f\cup g(a_1,\cdots, a_{p+q})=(-1)^{pq}
f(a_1,\cdots, a_p)\otimes_Ag(a_{p+1},\cdots a_{p+q}).$$
 
La formule 
$b(f\cup g)=bf\cup g+(-1)^pf\cup bg$ montre que $f\cup g$ est un $b$-cocycle d\`es que $f$ et $g$ le sont, ce qui permet d\'efinir $\alpha\cup \beta$ comme la classe de cohomologie du $b$-cocycle
$f\cup g$. 
 
\medskip
 
{\sl Le cap-produit.}
 
Il s'agit d'une application  $k$-lin\'eaire
$$\cap: H_r(A,N)\otimes_kH^p(A,M)\to H_{r-p}(A,N\otimes_AM).$$
 
\'Ecrivons $z\in H_r(A,N)$ sous la forme 
$z=cl(\overline z)$ où 
$\overline z=(n, a_1,\cdots, a_r)$ est un $b$-cycle de $C_r(A,N)$. Ecrivons  $\alpha\in H^p(A,M)$ sous la forme $\alpha=cl(f)$ où $f$  est un $b$-cocycle de  $C^p(A,M)$.
D\'efinissons  l'\'el\'ement $\overline z\cap f$ de 
$C_{r-p}(A, N\otimes_AM)$ par la relation 
$$\overline z\cap f=
(-1)^{rp}(n\otimes_Af(a_1,\cdots, a_p),a_{p+1},\cdots, a_r).$$
 
La relation 
$b(\overline z)\cap f=\overline z\cap bf+(-1)^p b(\overline z\cap f)$ montre que $\overline z\cap f$ est un $b$-cycle de \hbox{$C_{r-p}(A, N\otimes_AM)$.} Ceci permet de d\'efinir 
$z\cap \alpha=cl(\overline z\cap f)$ comme  la classe d'homologie de ce cycle. 

\medskip
Le {\sl cup}-produit et le {\sl cap}-produit sont reli\'es 
par la formule suivante, de d\'emonstration imm\'ediate \`a partir des d\'efinitions ci-dessus.
 
Pour $z\in H_r(A,N)$, $\alpha\in H^p(A,M)$ et 
$\beta \in H^q(A, M')$, dans 
$H_{r-(p+q)}(A,\ N\otimes_AM\otimes_AM')$, on a l'\'egalit\'e 
 
$$(z\cap\alpha)\cap\beta=z\cap(\alpha\cup\beta).$$
\medskip

{\sl Le crochet de Gerstenhaber.}

Soient $\alpha\in H^p(A,A)$ et $\beta\in H^q(A,A)$. 
\'Ecrivons $\alpha=cl(f)$ et $\beta=cl(g)$ avec $f$ et $g$ cocycles respectifs de $C^p(A,A)$  et de $C^q(A,A)$. Pour d\'efinir le crochet de Gerstenhaber $[\alpha,\beta] \in H^{p+q-1}(A,A)$ de $\alpha$ et $\beta$, on introduit les applications $k-$lin\'eaires
$f\circ_i g:A^{\otimes p+q-1}\to A$, d\'efinies pour $1\leq i\leq p$ par 
$$f\circ_i g(a_1,\cdots, a_{p+q-1})
=f(a_1,\cdots, a_{i-1},g(a_i,\cdots,a_{i+q-1}),
a_{i+q},\cdots a_{p-1+q}). $$ On pose ensuite
$$f\circ g=\sum_{i=1}^p(-1)^{(i-1)(q-1)}f\circ_ig$$  et enfin 
$$[f,g]=f\circ g-(-1)^{(p-1)(q-1)}g\circ f.$$
Gerstenhaber a montr\'e la formule 
$$b(f\circ g)=f\circ b(g)+(-1)^{q-1}b(f)\circ g
+(-1)^{q-1}(g\cup f-(-1)^{pq}f\cup g).$$
Ceci montre que le crochet $[f,g]\in C^{p+q-1}(A,A)$ de deux cocycles $f$ et $g$ est \'egalement un cocycle. La classe de cohomologie du cocycle $[f,g]$ est par d\'efinition le crochet $[\alpha,\beta]$ des classes $\alpha$ et $\beta$.
\medskip

M. Gerstenhaber a montr\'e
\medskip

{\bf Th\'eor\`eme 2.1.} ([Ge])
{\sl L'alg\`ebre de cohomologie de Hochschild $H^*(A,A)$ est une alg\`ebre de Gerstenhaber.}
\medskip

{\sl Le bord de Connes.} ([C], [L])

Le bord $B: C_r(A,A)\to C_{r-1}(A, A)$ est donn\'e par la formule 
$$B(a_0,a_1\cdots,a_r)=\sum_{j=0}^{r}(-1)^{jr}(1,a_j,\cdots, a_r,a_0,\cdots, a_{j-1}).$$
\medskip
Compte tenu de la relation $Bb+bB=0$, le bord de Connes induit un morphisme de $k$-espaces vectoriels, qu'on note \'egalement $B$ par abus de language :

$$B:H_r(A,A)\to H_{r+1}(A,A).$$

D'apr\`es un r\'esultat de D. Tamarkin et B. Tsygan, on a 
\medskip

{\bf Th\'eor\`eme  2.2.} ([T-T])
{\sl  Le triplet  $(H^*(A,A), H_*(A,A),B)$ est un calcul de Tamarkin-Tsygan.} 

\bigskip

{\bf 3. Alg\`ebres  de type $FP$.}
 \medskip
 
{\bf D\'efinition 3.1.}
{\sl Soient $k$ un corps et $A$ un $k$-alg\`ebre associative.  On dit que $A$ est une alg\`ebre de type $FP$ si l'alg\`ebre $A$ admet une r\'esolution projective de longueur finie par des $A^e$-modules projectifs de type fini.}
 \medskip
 
Pour une alg\`ebre $A$ de type $FP$, la dimension cohomologique de $A$ est 
$$d=\hbox{pdim}_{A^e}(A).$$
 On pose  
$${\cal D}=H^d(A, A^e).$$
 
\medskip
Rappelons ([B-T]) que le $k$-espace vectoriel ${\cal D}=H^d(A, A^e)$ est un $A^e$-module \`a gauche.
\medskip

On v\'erifie sans difficult\'e la proposition suivante.
\medskip

{\bf Proposition 3.2.}
{\sl Le cap-produit
$$H_d(A,{\cal D})\otimes_k\ H^d(A, A^e)\to
 H_0(A,\ {\cal D}\otimes_A A^e)$$
 $$z\otimes\alpha\to z\cap\alpha$$ fournit    un morphisme $k$-lin\'eaire 
$$H_d(A,{\cal D})\to Hom_{A^e}({\cal D},{\cal D}).$$
$$z\to z\cap -$$
}

\medskip

{\bf Proposition 3.3.}
{\sl Soit $A$ une alg\`ebre de type $FP$ de dimension cohomologique $d$. Pour tout $A^e$-module $M$, l'application
$$H_d(A,M)\to Hom_{A^e}({\cal D}, M)$$
$$z\to (z\cap-)_{\mid_{\cal D}} $$
 est un isomorphisme de $k$-espace vectoriels. 
}

\medskip

D\'emonstration. Le {\sl cap}-produit 
$$H_d(A,M)\otimes_kH^d(A, A^e)\to 
H_0(A, M\otimes_AA^e)\cong M$$
fournit un morphisme 
$$H_d(A,M)\to \hbox{Hom}({\cal D},M)$$ 
$$z\to (z\cap-)_{\mid_{\cal D}}$$
et on v\'erifie facilement que l'application $z\cap\ -$ est $A^e$-lin\'eaire.

\smallskip

Soit $P_*$ une $A^e$-r\'esolution projective de type fini de $A$, de longueur $d$. On a 
$$H_i(A,M)=H_i(P_*\otimes _{A^e}M).$$
En particulier on a la suite exacte courte 
$$0\to H_d(A,M)\to P_d \otimes_{A^e}M\to P_{d-1}\otimes_{A^e}M.\leqno (1)$$

\smallskip

Posons 
$\overline P^*=Hom_{A^e}(P_*,A^e)$. On a $H^i(\overline P^*)=H^i(A, A^e).$
En particulier, on a la suite exacte courte 
$$\overline P^{d-1}\to 
\overline P^d\to H^d(A, A^e)={\cal D}\to 0.\leqno (2)$$

\smallskip

Par application du foncteur 
$\hbox{Hom}_{A^e}(-,M)$ \`a la suite exacte courte $(2)$, on obtient la suite exacte 
$$0 \to \hbox{Hom}_{A^e}({\cal D}, M)\to 
\hbox{Hom}_{A^e}(\overline P^d, M)\to
\hbox{Hom}_{A^e}(\overline P^{d-1}, M).\leqno (3)$$

\smallskip

Puisque les $P_i$ sont projectifs de type finis, on a  des isomorphismes 
$$\hbox{Hom}_{A^e}(\overline P^i, M)\cong P_i\otimes_{A^e}M$$ et la suite exacte $(3)$ s'\'ecrit donc 
$$0\to \hbox{Hom}_{A^e}({\cal D},M)\to
P_d\otimes _{A^e}M\to P_{d-1}\otimes_{A^e}M.\leqno (4)$$
\smallskip

Les suites exactes $(1)$ et $(4)$ fournissent l'isomorphisme 
$H_d(A,M)\cong Hom_{A^e}({\cal D}, M).$

\medskip

Remarque : Pour $M={\cal D}$, la proposition 3.3  montre que  
$$H_d(A,{\cal D})\to Hom_{A^e}({\cal D},{\cal D}).$$
$$z\to z\cap\ -$$
est un isomorphisme. 
Ceci conduit \`a la d\'efinition suivante.

\medskip

{\bf D\'efinition 3.4.}
{\sl Soit $A$ une alg\`ebre de type $FP$, de dimension cohomologique $d$. L'unique \'el\'ement $c$ de $H_d(A,{\cal D})$ tel que 
$$(c\cap-)_{\mid_{\cal D}} =\hbox{id}_{\cal D}$$ s'appelle la classe fondamentale de l'alg\`ebre  $A$.}
\medskip

{\bf Proposition 3.5.}
{\sl Soit $A$ une alg\`ebre de type $FP$, de dimension cohomologique $d$, de classe fondamentale $c\in H_d(A, {\cal D})$.
Pour tout $A^e$-module $M$, le cap-produit 
$$c\cap-: H^d(A;M)\to H_0(A,{\cal D}\otimes_AM)$$ est un isomorphisme.}
\medskip

D\'emonstration. Si $M=A^e$ est libre de rang $1$, par d\'efinition de la classe fondamentale
le {\sl cap}-produit $c\cap-$ est l'application 
$\hbox{id}: {\cal D}\to {\cal D}$.
\smallskip

Pour traiter le cas des modules libres, on regarde $c\cap-$ comme une transformation naturelle du foncteur $H^d(A,-)$ vers le foncteur $H_0(A,{\cal D}\otimes_A\ -)$. Ces deux foncteurs sont additifs. En outre, comme tout foncteur $Tor$, le foncteur 
$H_0(A,{\cal D}\otimes_A-)$ commute aux limites directes. Puisque $A$ est de type $FP$, le foncteur 
$\hbox{Ext}^d_{A^e}(A, -)$ commute \'egalement aux limites directes ([Br], VIII, 4.8, p. 196). Ceci montre que si $M$ est un $A^e$-module libre, le {\sl cap}-produit $c\cap-: H^d(A,M)\to H_0(A,{\cal D}\otimes_AM)$ est un isomorphisme.
\smallskip

Enfin, si $M$ est un $A^e$-module quelconque, on obtient le r\'esultat gr\^ace \`a une suite exacte $F'\to F\to M\to 0$, où $F$ et $F'$ sont libres.

\vfill\eject

 {\bf 4. Dualit\'e de Van den Bergh.}
\medskip

{\bf D\'efinition 4.1.}
{\sl Soit $A$ une alg\`ebre de type $FP$ de dimension cohomologique $d$. On dit que $A$ est une alg\`ebre \`a dualit\'e de Van den Bergh si 
$H^i(A, A^e)=0$ pour $i\not= d$
et si  
${\cal D}=H^d(A, A^e)$ est un $A^e$-module inversible.

\smallskip

Dans ce cas, $\cal D$ est appel\'e le module dualisant de l'alg\`ebre \`a dualit\'e de Van den Bergh $A$.}

\medskip

{\bf Th\'eor\`eme 4.2.}
{\sl Soit  $A$ une alg\`ebre \`a dualit\'e de Van den Bergh, de module dualisant ${\cal D}=H^d(A, A^e)$, de classe fondamentale 
$c\in H_d(A,{\cal D})$. Alors, pour tout $A^e$-module $M$ et pour tout entier $p\geq 0$, le cap-produit 
$$c\cap- :H^p(A, M)\to H_{d-p}(A, {\cal D}\otimes_AM)$$
est un isomorphisme.}
\medskip

D\'emonstration.
Introduisons les foncteurs 
$H_i=H_i(A,{\cal D}\otimes_A -)$ et 
$T_i=H^{d-i}(A, -)$.
De mani\`ere \'evidente, le foncteur $T_i$ est un foncteur homologique. 
\smallskip

{\sl Premi\`ere \'etape} : Grâce \`a $H^i(A, A^e)=0$ pour $i\not=d$, montrons que le foncteur $T_i$ est  
effa\c cable pour $i>0$ ([Br], III, 6, p. 72). Par hypoth\`ese, on a 
$T_i(A^e)=0$ pour $i>0$. Par additivit\'e du foncteur $T_i$, on en d\'eduit 
$T_i(L)=0$ pour $i>0$ et $L$ libre de rang fini, donc 
$T_i(P)=0$, pour $i>0$ et $P$ projectif de type fini.
D'apr\`es [Br], VIII, 4.6, p. 195, $T_i$ commute aux limites directes.  On en d\'eduit 
$T_i(L)=0$, pour $i>0$ et $L$ libre et donc aussi
$T_i(P)=0$ pour $i>0$ et $P$ projectif.

\smallskip

{\sl Deuxi\`eme \'etape} : Gr\^ace \`a $\cal D$ inversible, montrons \`a pr\'esent que 
$H_i$ est un foncteur homologique. 
Soit $$0\to M'\to M\to M''\to 0$$ une suite exacte courte de $A^e$-modules. Puisque $\cal D$ est inversible, $\cal D$ est un $A$-module projectif \`a droite, donc plat. On  a par cons\'equent la suite exacte courte de $A^e$-modules
$$0\to{\cal D}\otimes_AM'\to
{\cal D}\otimes_AM\to
{\cal D}\otimes_AM''\to 0.$$ 
Cette suite exacte courte fournit la suite exacte longue 
\[\cdots\to H_{i+1}(A,{\cal D}\otimes_AM'')\to 
H_i(A,{\cal D}\otimes_AM')\to
H_i(A,{\cal D}\otimes_AM)\to
H_i(A,{\cal D}\otimes_AM'')\to
\cdots\]
c'est-\`a-dire qu'on a la suite exacte longue
$$¤\cdots\to H_{i+1}(M'')\to 
H_i(M')\to H_i(M)\to H_i(M'')
\to H_{i-1}(M')\to \cdots$$
ce qui montre que $H_i$ est un foncteur homologique.
\smallskip

{\sl Troisi\`eme \'etape} : Montrons \`a pr\'esent que $H_i$ est un foncteur effa\c cable pour $i>0$. Puisque $\cal D$ est $A^e$-inversible, il est $A$-projectif \`a gauche. Donc si $P$ est un $A^e$-module projectif, ${\cal D}\otimes_AP$ est encore  projectif et par cons\'equent 
$H_i(P)=H_i(A,{\cal D}\otimes_AP)$ est nul  pour $i>0$.
\smallskip

{\sl Quatri\`eme \'etape} : D'apr\`es la proposition 3.5, $c\cap- : T_0\to H_0$ est un isomorphisme.

\smallskip

{\sl Cinqui\`eme \'etape} : 
Par d\'ecalage d'indice (cf [Br], III, 7.3, p. 75), on d\'eduit de tout ce qui pr\'ec\`ede que pour tout $i>0$, $c\cap -: T_i\to H_i$ est un isomorphisme. 
\bigskip

\bigskip

{\bf 5. Une d\'emonstration du th\'eor\`eme de V. Ginzburg.}
\medskip

{\bf D\'efinition 5.1.}
{\sl On dit que l'alg\`ebre $A$ est une alg\`ebre de Calabi-Yau  de dimension $d$ si 
$A$ est une alg\`ebre \`a dualit\'e de Van den Bergh de dimension cohomologique $d$ dont le module dualisant 
${\cal D}=H^d(A, A^e)$ est un $A^e$-module isomorphe \`a $A$.}
\medskip

 Soit $A$ une alg\`ebre de Calabi-Yau de dimension  $d$. Posons \hbox{$H^*=H^*(A,A)$} et 
$H_*=H_*(A, A)$. 
Soit 
$B$  le bord  de Connes. D'apr\`es le r\'esultat de Tamarkin et Tsygan rappel\'e en 2.2,   $(H^*,H_*,\  B)$ est un calcul de Tamarkin-Tsygan.
Soit $c\in H_d(A,{\cal D})$ la  classe fondamentale de l'alg\`ebre $A$. Du th\'eor\`eme 4.2, nous d\'eduisons que $(H^*,H_*,\  B,\ c)$ est un calcul de Tamarkin-Tsygan \`a dualit\'e. Compte tenu du  th\'eor\`eme 1.6, on a donc montr\'e
\medskip

{\bf Th\'eor\`eme 5.2.} (V. Ginzburg)
{\sl Une alg\`ebre de Calabi-Yau   est une BV-alg\`ebre. Plus pr\'ecis\'ement,
soit $A$ une alg\`ebre de Calabi-Yau de dimension $d$, de classe fondamentale  $c$. Soit 
$D=c\cap\ -\ $ l'isomorphisme de dualit\'e de Van den Bergh.  Alors $\Delta:= (-1)^d\ DB D^{-1}$ est un g\'en\'erateur du crochet de Gerstenhaber de $H^*(A,A)$, c'est-\`a-dire que pour tout 
$\alpha\in H^p(A,A)$ et tout $\beta\in H^q(A,A)$, on a l'\'egalit\'e 
$$[\alpha,\beta]=
\Delta(\alpha\cup\beta)-(-1)^p\alpha\cup\Delta(\beta)-\Delta (\alpha)\cup\beta.$$}

\bigskip
 
{\bf 6. Une d\'emonstration d'un r\'esultat de M. Kontsevich.}
\medskip

V. Ginzburg attribue le r\'esultat suivant \`a M. Kontsevich.
\medskip
 
{\bf Th\'eor\`eme 6.1.} ([G], 6.1.1).
{\sl Soit $X$ une vari\'et\'e orient\'ee asph\'erique de dimension $3$. Alors ${\bf C}[\pi_1(X)]$, alg\`ebre du groupe fondamental de $X$ est une alg\`ebre de Calabi-Yau de dimension $3$.}
\medskip

Nous proposons une d\'emonstration de ce r\'esultat sous la forme suivante.
\medskip
   
{\bf Lemme 6.2.}
{\sl Soit $G$ un groupe \`a dualit\'e de Poincar\'e. On suppose $G$ orientable de dimension cohomologique $d$. Soit $k$ un corps de caract\'eristique z\'ero ou premi\`ere \`a l'ordre du groupe $G$. Alors l'alg\`ebre du groupe $A=k[G]$ est une alg\`ebre de Calabi-Yau de dimension $d$.}
\medskip
 
Pour d\'emontrer ce lemme, rappelons le vocabulaire des groupes \`a dualit\'e de Poincar\'e ([B-E], [Br],VIII,§ 10). 
\medskip

On dit qu'un groupe $G$ est de type $FP$ si ${\mathbb Z}$ admet une 
${\mathbb Z}[G]$-r\'esolution de longueur finie par des ${\mathbb Z}[G]$-modules projectifs de type fini. Dans ce cas, la dimension cohomologique $d$ du groupe $G$ est  \hbox{$d=\hbox{pdim}_{{\mathbb Z}[G]}{\mathbb Z}.$}
\medskip

Le groupe $G$ de type $FP$, de dimension cohomologique $d$  est \`a dualit\'e de Poincar\'e si
 $\hbox{Ext}_{{\mathbb Z}[G]}^i({\mathbb Z},{\mathbb Z}[G])=0$ pour $i\not=0$ et si 
 ${\cal D}_G:=\hbox{Ext}_{{\mathbb Z}[G]}^d({\mathbb Z},{\mathbb Z}[G])$ est un $\mathbb Z$-module isomorphe \`a $\mathbb Z$. Le ${\mathbb Z}[G]$-module ${\cal D}_G$ s'appelle le module dualisant du groupe $G$. 
 
\medskip
  
Le groupe $G$  est \`a dualit\'e de Poincar\'e orientable s'il est \`a dualit\'e de Poincar\'e et si  son module dualisant ${\cal D}_G$ est un ${\mathbb Z}[G]$-module isomorphe au 
 ${\mathbb Z}[G]$-module trivial  $\mathbb Z$.
 \medskip

 \medskip

{\bf Proposition 6.3.}
 {\sl Soient  $G$ un groupe,  $k$ un corps de caract\'eristique z\'ero ou premi\`ere \`a l'ordre de $G$, et soit $A$ l'alg\`ebre du groupe $k[G]$. On suppose que  $G$ est un groupe de type $FP$ de dimension cohomologique $d$. Alors l'alg\`ebre $A$ est une alg\`ebre de type $FP$ de dimension cohomologique $d$ et pour tout $i$ on a des isomorphismes de $A^e$-modules
 $$\hbox{Ext}_{A^e}^i(A,A^e)
 \cong 
\hbox{Ext}^i_{{\mathbb Z}[G]}({\mathbb Z},{\mathbb Z}[G])\otimes_{{\mathbb Z}}A.$$
 
\smallskip
 
 En particulier, si $G$ est \`a dualit\'e de Poincar\'e de module dualisant ${\cal D}_G$, alors $A$ est une alg\`ebre \`a dualit\'e de Van den Bergh dont le module dualisant
 ${\cal D}_A$ est $A^e$-isomorphe \`a ${\cal D}_G\otimes_{\mathbb Z}A$.} 
 
\medskip

D\'emonstration. 
Puisque $k$ est de caract\'eristique z\'ero ou premi\`ere \`a l'ordre de $G$, $A$ est un 
${\mathbb Z}[G]$-module plat. Le foncteur 
$$-\ \otimes_{{\mathbb Z}[G]}A:{\mathbb Z}[G]\hbox{-Mod}\to A\hbox{-Mod}$$
est donc exact. Puisque $A$ est un $k$-module projectif, il est plat. Le foncteur 
$$-\ \otimes_kA:A\hbox{-Mod}\to A^{\otimes 2}\hbox{-Mod}$$ est \'egalement exact.
 Le foncteur compos\'e 
 $-\ \otimes_{{\mathbb Z}[G]}A^{\otimes 2} $ est donc exact. Puisque $G$ est de type FP, de dimension $d$, il existe une 
 ${\mathbb Z}[G]$-r\'esolution de longueur $d$
 $$0\to P_d\to\cdots\to P_0\to^\varepsilon{\mathbb Z}\to 0$$
 du ${\mathbb Z}[G]$-module trivial $\mathbb Z$ par des ${\mathbb Z}$-modules projectifs de type fini. Par application du foncteur exact 
 $-\ \otimes_{{\mathbb Z}[G]}A^{\otimes 2}$ et compte tenu de l'isomorphisme de $k$-alg\`ebres $A^e\cong A^{\otimes 2}$, on obtient une $A^e$-r\'esolution de longueur $d$ de $A$ par des $A^e$-modules projectifs de type fini. Ceci  montre que $A$ est de type FP, de dimension cohomologique $d$.
 
Les isomorphismes
$$
\hbox{Ext}^i_{A^e}(A, A^e)
\cong 
\hbox{Ext}^i_{{\mathbb Z}[G]}({\mathbb Z},{\mathbb Z}[G])\otimes_{{\mathbb Z}[G]}A^e
\cong
\hbox{Ext}^i_{{\mathbb Z}[G]}({\mathbb Z},{\mathbb Z}[G])\otimes_{{\mathbb Z}[G]}A^{\otimes 2}
\cong
\hbox{Ext}^i_{{\mathbb Z}[G]}({\mathbb Z},{\mathbb Z}[G])\otimes_{\mathbb Z}A$$
ach\`event la d\'emonstration de la proposition 6.3.
\medskip

D\'emonstration du lemme 6.2. D'apr\`es  6.3, puisque $G$ est \`a dualit\'e de Poincar\'e, $A=k[G]$ est \`a dualit\'e de Van den Bergh et
$${\cal D}_A
\cong 
{\cal D}_G\otimes_{\mathbb Z}A
.$$ Puisque $G$ est orientable, l'action de $G$ sur ${\cal D}_G=\mathbb Z$ est triviale et  on a donc un isomorphime de $A^e$-modules 
${\cal D}_A\cong A$,
ce qui montre que $A$ est Calabi-Yau.
\medskip

{\bf Remarques.} 
\smallskip

- De 6.2 et 5.1, on d\'eduit que sous les hypoth\`eses 6.2, l'alg\`ebre de cohomologie de Hochschild de $A=k[G]$ est une BV-alg\`ebre.  Dans le cas des vari\'et\'es orient\'ees asph\'eriques de dimension $d$, D. Vaintrob [V] a montr\'e le lien entre cette structure BV sur $H^*(A,A)$ et la structure BV de Chas-Sullivan sur l'homologie singuli\`ere $H_*({\cal L}(X))$ de l'espace des lacets libres sur $X$.
\smallskip

- Si $G$ est \`a dualit\'e de Poincar\'e {\sl non orientable}, l'alg\`ebre $A=k[G]$ est toujours une alg\`ebre \`a dualit\'e de Van den Bergh mais puisque  l'op\'eration de $G$ sur ${\cal D}_G={\mathbb Z}$ {\sl n'est pas} triviale, le module dualisant ${\cal D}_A$ de l'alg\`ebre $A$  est un $A^e$-module isomorphe  au module tordu  ${}_\varphi A$, (c'est-\`a-dire $g\cdot x\cdot h=\varphi(g)xh$)  où 
$\varphi $ est l'isomorphime $\varphi: A\to A$
d\'efini par $\varphi (g)=n_g g$ avec $n_g=g\cdot 1$.

\bigskip
\bigskip

{\bf Bibliographie.}
\bigskip

[B-T],
R. Berger \& R. Taillefer, 
Poincar\'e-Birkhoff-Witt deformations of Calabi-Yau algebras,
{\sl J. Noncommutative Geom.}, {\bf 1}, 2007, 241-270.
\medskip

[B-E],
R. Bieri \& B. Eckmann,
Groups with homological duality generalizing Poincar\'e duality,
{\sl Inventiones Math.}, {\bf 20}, 1973, 103-124.
\medskip

[Bi],
N. Bourbaki,
{\sl Alg\`ebre}, chap. 10, Alg\`ebre homologique, Masson, 1980.
\medskip

[Br],
K. Brown, 
{\sl Cohomology of groups}, GTM 87, 1982, Springer.
\medskip

[C],
A. Connes,
Non-Commutative differential geometry, 
{\sl Publ. Math. IHES}, {\bf 62}, 1985, 257-360.
\medskip

[Ge], 
M. Gerstenhaber,
The cohomology structure of an associative ring,
{\sl Ann. of Math.}, {\bf 78}, 1963, 267-288.
\medskip

[Gi],
V. Ginzburg,
Calabi-Yau algebras, arXiv: 06 12 139, 2006.

\medskip

[H],
J. Huebschmann,
Duality for Lie-Rinehart algebras and the modular class,
{\sl J. Reine Angew. Math.}, {\bf 510}, 1999, 103-159.
\medskip

[L],
J.-L. Loday,
{\sl Cyclic homology},
Grundlehren der Math. Wissenschaften, 301, Springer, 1992.
\medskip

[T-T],
D. Tamarkin \& B. Tsygan,
{The ring of differential operators on forms in non-commutative calculus}, {\sl Proceedings of Symposia in Pure Math., Am. Math. Soc.}, vol. 73, 2005, 105-131. \medskip

[V],
D. Vaintrob,
The string topology BV-algebra, Hochschild cohomology and the Goldmann bracket on surfaces, Arxiv, 07 02 859, 2007.

\medskip

[VdB],
M. Van den Bergh,
A relation between Hochschild homology and cohomology for Gorenstein rings,
{\sl Proceedings of the Am. Math. Soc.}, {\bf 126}, 1998, 1345-1348 and
  {\bf 130}, 2002,  2809-2810 (Erratum).
\medskip

[X], P. Xu,
Gerstenhaber algebras and BV-algebras in Poisson geometry,
{\sl Com. Math. Phys.}, {\bf 200}, 1999, 545-560. 
\bigskip

\end{document}